\input amstex
\documentstyle{amsppt}
\magnification=1200
\voffset-1cm
\topmatter
\title On Vassiliev Invariants for Links in Handlebodies
\endtitle
\author V. V. Vershinin \endauthor
\address
Sobolev Institute of Mathematics, Novosibirsk, 630090, Russia
\endaddress
\email versh\@math.nsc.ru \endemail
\thanks The research was supported in part by RFFR Grants 96-01-016331 and
96-01-00062G
\endthanks
\abstract The notion of Vassiliev algebra in case of hanlebodies is
developed. The analogues of the results of John Baez for links in
handlebodies are proved. That means that there exists a one-to-one
correspondence between the special class of finite type invariants of links
in hanlebodies and the homogeneous Markov traces on Vassiliev algebras.
This approach uses the singular braid monoid and braid group in a handlebody
and the generalizations of the theorem of J.~Alexander and the theorem of
A.~A.~Markov for singular links and braids and the relative version of
Markov's theorem.
\endabstract
\endtopmatter
\document
\specialhead \S 1. Introduction
\endspecialhead
Knot and link invariants introduced by V.~A.~Vassiliev [1] attracted a lot
of attention and are studied deeply from various points of view as well as
their connections with other fields of mathematics. In spite of the short
period of time there exists a big amount of works on this subject.

The conception of finite type invariants was applied also for study of knots
and links in 3-manifolds. Here Xiao-Song Lin [2] considered contractible
manifolds and homological spheres and proved that the module of Vassiliev
invariants of the order $n$ for such a manifold coincides wih that of the
standard 3-dimensional euclidean space. V.~A.~Vassiliev continues his
original approach considering the space of smooth maps from a circle to a
given manifold, its discriminant, cohomology and spectral sequence [3].
V.~Gorynov in his work [4] in particular develops the formalism of chord
diagrams for solid torus and proves the analogue of M.~Kontsevich theorem in
this case.

The classical theorems of J.~Alexander and A.~A.~Markov established the close
interconnections between the study of links and the theory of braids. For
Vassiliev invariants relations between singular braids and singular knots
arise. They were studied by John Baez [5] and Joan Birman [6].

The aim of this paper is to carry the results of John Baez [5] for links
in handlebodies. In our approach we use the braid group in a handlebody,
which was studied by A.~B.~Sossinsky [7] and the author [8, 9]. The
essential ingredients are the generalizations of the theorem of
J.~Alexander and the theorem of A.~A.~Markov for singular links and braids
proved by Joan Birman [6] and Bernd Gemein [10] and the relative version of
Markov's theorem proved by Sofia Lambropoulou and Colin Rourke [11]. This
version allows us to have Markov's theorem for links in hanlebodies. The
proofs of John Baez [5] remain valid after necessary adaptations to the case
of a handlebody.

The author is thankful to Sofia Lambropoulou, Bernd Gemein and Viktor
Vassiliev for supplying him with the preprints of their papers.

\specialhead \S 2. Singular braids in handlebodies
\endspecialhead
Let $SB_n$ be the singular braid monoid on $n$ strings [5, 6], which is
also called the Baez-Birman monoid or generalized braid monoid. It has the
generators
$\bar\sigma_i, \ \bar\sigma_i^{-1}$, $\bar a_i$,  $i=1, ..., n-1,$ and
relations
$$\align &\bar\sigma_i\bar\sigma_j=\bar\sigma_j\bar\sigma_i, \ \text {if} \
|i-j| >1,\\
&\bar a_i\bar a_j=\bar a_j\bar a_i, \ \text {if} \  |i-j| >1,\\
&\bar a_i\bar\sigma_j=\bar\sigma_j\bar a_i, \ \text {if} \ |i-j| \not=1,\\
&\bar\sigma_i\bar\sigma_{i+1}\bar\sigma_i=\bar\sigma_{i+1}\bar\sigma_i
\bar\sigma_{i+1}, \tag1 \\
&\bar\sigma_i\bar\sigma_{i+1}\bar a_i=
\bar a_{i+1}\bar\sigma_i\bar\sigma_{i+1},\\
&\bar\sigma_{i+1}\bar\sigma_i\bar a_{i+1}=
\bar a_i\bar\sigma_{i+1}\bar\sigma_i,\\
&\bar\sigma_i\bar\sigma_i^{-1}=\bar\sigma_i^{-1}\bar\sigma_i=1. \endalign $$
Geometrically $\bar\sigma_i$ corresponds to the canonical generator of the
braid group (right-handed crossing) and $\bar a_i$ represents an intersection
of the $i$th and $(i+1)$st strand, just as in Figures 1 and 2. More detailed
geometric interpretation of the Baez-Birman monoid can be found in the paper
of Joan Birman [6].

\input lpic
\noindent\centerline{\latexpic(0,130)(0,-15)
\thicklines
\put(-100,100){\lline(0,-1){100}}
\put(-50,100){\lline(0,-1){100}}
\put(-25,100){\lline(1,-2){50}}
\put(25,100){\lline(-1,-2){20}}
\put(-25,0){\lline(1,2){20}}
\put(50,100){\lline(0,-1){100}}
\put(100,100){\lline(0,-1){100}}
\put(-100,110){\makebox(0,0)[cc]{$1$}}
\put(-50,110){\makebox(0,0)[cc]{$i-1$}}
\put(-25,110){\makebox(0,0)[cc]{$i$}}
\put(25,110){\makebox(0,0)[cc]{$i+1$}}
\put(50,110){\makebox(0,0)[cc]{$i+2$}}
\put(100,110){\makebox(0,0)[cc]{$n$}}
\put(-75,50){\makebox(0,0)[cc]{.\quad.\quad.}}
\put(75,50){\makebox(0,0)[cc]{.\quad.\quad.}}
\endlatexpic}
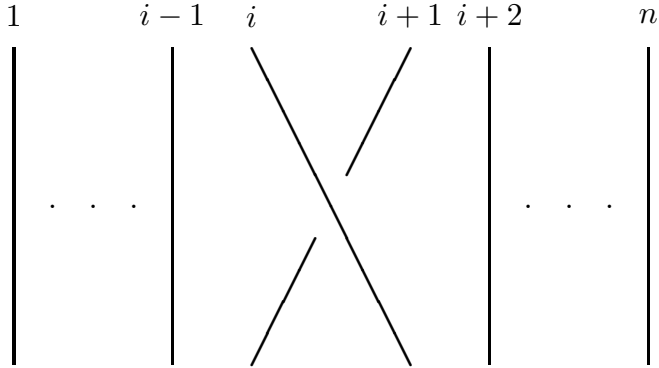
\centerline {Fig. 1}
\vglue 0.5cm
\penalty-10000

\noindent\centerline{\latexpic(0,130)(0,-10)
\thicklines
\put(0,50){\circle*{5}}
\put(-100,100){\lline(0,-1){100}}
\put(-50,100){\lline(0,-1){100}}
\put(-25,100){\lline(1,-2){50}}
\put(25,100){\lline(-1,-2){50}}
\put(50,100){\lline(0,-1){100}}
\put(100,100){\lline(0,-1){100}}
\put(-100,110){\makebox(0,0)[cc]{$1$}}
\put(-50,110){\makebox(0,0)[cc]{$i-1$}}
\put(-25,110){\makebox(0,0)[cc]{$i$}}
\put(25,110){\makebox(0,0)[cc]{$i+1$}}
\put(50,110){\makebox(0,0)[cc]{$i+2$}}
\put(100,110){\makebox(0,0)[cc]{$n$}}
\put(-75,50){\makebox(0,0)[cc]{.\quad.\quad.}}
\put(75,50){\makebox(0,0)[cc]{.\quad.\quad.}}
\endlatexpic}
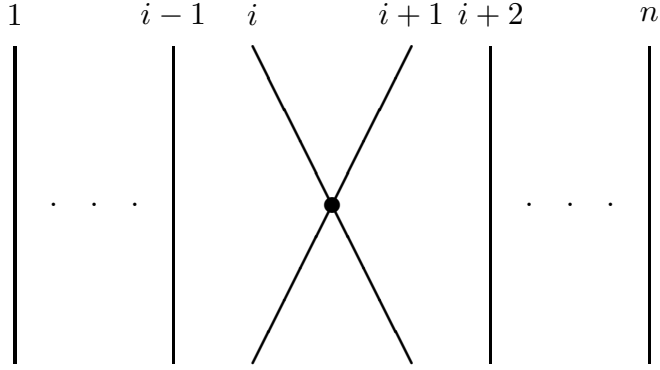
\centerline {Fig. 2}
\vglue 0.5cm

As the analogous object in a handlebody of the genus $g$ it is natural to
consider the submonoid of $SB_{g+n}$ in which the first $g$ strings are
unbraided and noncrossing. We denote it by $SB_n^g$. The braid group $Br_n^g$
in a handlebody of the genus $g$ was studied by A.~B.~Sossinsky [7] and the
author [8, 9]. It has the generators
$\tau_k, \ k=1,2,..., g$ and $\sigma_{1}, ..., \sigma_{n-1}$ and relations
$$\cases\sigma_i\sigma_j &=\sigma_j\sigma_i \ \ \text{if} \ \ |i-j| >1,\\
\sigma_i \sigma_{i+1} \sigma_i &= \sigma_{i+1} \sigma_i \sigma_{i+1},\\
\tau_k\sigma_i &=\sigma_i\tau_k \ \ \text{if} \ \ k\geq 1, \ i\geq 2,\\
\tau_k\sigma_1\tau_k\sigma_1&=\sigma_1\tau_k\sigma_1\tau_k, \ k=1,2,..., g;\\
 \tau_k\sigma_1^{-1}\tau_{k+l}\sigma_1 &=
\sigma_1^{-1}\tau_{k+l}\sigma_1\tau_k,
\ k=1,2,..., g-1; \ l=1,2,..., g-k. \\
\endcases \tag2 $$
The group $Br_n^g$ is considered as a subgroup of the classical braid
group $Br_{g+n}$ on $g+n$ strings, such that the braids from $Br_n^g$ leave
the first $g$ strings unbraided. We denote also by $\bar\sigma_j$ the
standard generators of the group $Br_{g+n}$. Then $\tau_k, \ k=1,2,..., g$,
are the following braids:
$$\tau_k = \bar\sigma_g\bar\sigma_{g-1}...\bar\sigma_{k+1}
\bar\sigma_{k}^2\bar\sigma_{k+1}^{-1}...\bar\sigma_{g-1}^{-1}
\bar\sigma_g^{-1}.$$
Geometrically such a braid (for $k=1$) is depicted in the Figure 3.

\vglue 0.3cm

\centerline{\noindent\latexpic(0,130)(-40,-15)
\thicklines
\put(-100,100){\lline(0,-1){57}}
\put(-100,0){\lline(0,1){30}}
\put(-80,100){\lline(0,-1){100}}
\put(-45,100){\lline(0,-1){100}}
\put(-25,100){\lline(-2,-1){15}}
\put(-85,70){\lline(-2,-1){10}}
\put(-50,88){\lline(-2,-1){25}}
\put(-50,12){\lline(-2,1){25}}
\put(-25,0){\lline(-2,1){15}}
\put(-5,100){\lline(0,-1){100}}
\put(30,100){\lline(0,-1){100}}
\put(-125,50){\lline(2,1){20}}
\put(-125,50){\lline(2,-1){40}}
\put(-100,110){\makebox(0,0)[cc]{$1$}}
\put(-80,110){\makebox(0,0)[cc]{$2$}}
\put(-45,110){\makebox(0,0)[cc]{$g$}}
\put(-25,110){\makebox(0,0)[cc]{$g+1$}}
\put(30,110){\makebox(0,0)[cc]{$g+n$}}
\put(-61,0){\makebox(0,0)[cc]{...}}
\put(13,0){\makebox(0,0)[cc]{...}}
\endlatexpic}
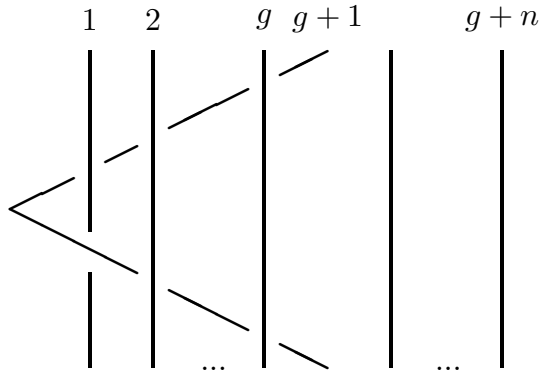
\centerline {Fig. 3}
\vglue 0.3cm
\penalty-10000
The elements $\sigma_i\in Br_n^g$ correspond to $\bar\sigma_{i+g}$.

The monoid $SB_n^g$ can be considered as a submonoid of $SB_{g+n}$ generated
by the elements $\sigma_{i}=\bar\sigma_{i+g}$, $\tau_k, \ k=1,2,..., g$ and
$a_{i}=\bar a_{i+g}$.

It is proved by R. Fenn, E. Keyman and C. Rourke [12] that the Baez-Birman
monoid $SB_k$ embeds in a group $SG_k$ which they call the {\it singular
braid group}:
$$SB_n\to SG_n.$$
It means that the elements $\bar a_i$ become invertible and all the relations
of $SB_k$ remain true.

It seems to be rather complicated to write a presentation of monoid $SB_n^g$,
but some relations can be obtained easily.
Let us introduce the following notations:
$$\tau_{i,m} =\tau_{i}\tau_{i+1}...\tau_{m}, \quad i\leq m.$$
\proclaim{Proposition 1} The following relations are fulfilled in the monoid
$SB_{g+n}$ and hence in the group $SG_{g+n}$:
$$\tau_{i,g}\bar\sigma_{g+1}\tau_{i,g}\bar a_{g+1}=\bar a_{g+1}\tau_{i,g}
\bar\sigma_{g+1}\tau_{i,g}, \quad i\leq g.$$
\endproclaim
\demo{Proof} These relations are obtained in the following way. We express
the elements $\bar a_i$, $i\leq g$, by means of $\bar a_{g+1}$ and
$\bar\sigma_j$ using firstly the fifth relation in (1) and secondly the sixth
relation in (1). Then we equate these two expressions. We start with
$\bar a_g$:
$$\bar a_g=\bar\sigma_{g+1}\bar\sigma_g\bar a_{g+1}
\bar\sigma_{g}^{-1}\bar\sigma_{g+1}^{-1},$$
$$\bar a_g=\bar\sigma_{g+1}^{-1}\bar\sigma_g^{-1}\bar a_{g+1}
\bar\sigma_{g}\bar\sigma_{g+1}.$$
Hence we have
$$\bar\sigma_g\bar\sigma_{g+1}^2\bar\sigma_g\bar a_{g+1}=
\bar a_{g+1}\bar\sigma_{g}\bar\sigma_{g+1}^2\bar\sigma_g.\tag3$$
Let us multiply on the right side both parts of the equality (3) by
$\bar\sigma_{g+1}$ and then use the braid group relations and the relation
$\bar\sigma _{g+1}\bar a_{g+1}=\bar a_{g+1}\bar\sigma _{g+1}$. We get the
following equalities for the left part of (3):
$$\bar\sigma_{g+1}\bar\sigma_g\bar\sigma_{g+1}^2\bar\sigma_g\bar a_{g+1}=
\bar\sigma_{g}\bar\sigma_{g+1}\bar\sigma_{g}\bar\sigma_{g+1}\bar\sigma_g
\bar a_{g+1}=\bar\sigma_g^2\bar\sigma_{g+1}\bar\sigma_g^2\bar a_{g+1},$$
and analogously for the right part:
$$\bar\sigma_{g+1}\bar a_{g+1}\bar\sigma_{g}\bar\sigma_{g+1}^2\bar\sigma_g=
\bar a_{g+1}\bar\sigma_g^2\bar\sigma_{g+1}\bar\sigma_g^2.$$
Using the expressions for $\tau_k$ we obtain:
$$\tau_{g,g}\bar\sigma_{g+1}\tau_{g,g}\bar a_{g+1}=\bar a_{g+1}\tau_{g,g}
\bar\sigma_{g+1}\tau_{g,g}.$$
In general case we have the following expressions
$$\bar a_1=\bar\sigma_2\bar\sigma_1\bar\sigma_3\bar\sigma_2...\bar\sigma_g
\bar a_{g+1}\bar\sigma_g^{-1}...\bar\sigma_2^{-1}\bar\sigma_3^{-1}
\bar\sigma_1^{-1}\bar\sigma_2^{-1},$$
$$\bar a_2=\bar\sigma_3\bar\sigma_2...\bar\sigma_g\bar a_{g+1}
\bar\sigma_g^{-1}...\bar\sigma_2^{-1}\bar\sigma_3^{-1},$$
$$...,$$
$$\bar a_i=\bar\sigma_{i+1}\bar\sigma_i\bar\sigma_{i+2}\bar\sigma_{i+1}...
\bar\sigma_g\bar a_{g+1}\bar\sigma_g^{-1}...\bar\sigma_{i+1}^{-1}
\bar\sigma_{i+2}^{-1}\bar\sigma_i^{-1}\bar\sigma_{i+1}^{-1}, \tag4 $$
$$...,$$
$$\bar a_{g-1}=\bar\sigma_{g}\bar\sigma_{g-1}\bar\sigma_{g+1}\bar\sigma_g
\bar a_{g+1}\bar\sigma_{g}^{-1}\bar\sigma_{g+1}^{-1}\bar\sigma_{g-1}^{-1}
\bar\sigma_{g}^{-1},$$
$$\bar a_1=\bar\sigma_2^{-1}\bar\sigma_1^{-1}\bar\sigma_3^{-1}
\bar\sigma_2^{-1}...\bar\sigma_g^{-1}\bar a_{g+1}\bar\sigma_g... \bar\sigma_2
\bar\sigma_3\bar\sigma_1\bar\sigma_2, $$
$$\bar a_2=\bar\sigma_3^{-1}\bar\sigma_2^{-1}...\bar\sigma_g^{-1}\bar a_{g+1}
\bar\sigma_g...\bar\sigma_2\bar\sigma_3,$$
$$...,$$
$$\bar a_i=\bar\sigma_{i+1}^{-1}\bar\sigma_i^{-1}\bar\sigma_{i+2}^{-1}
\bar\sigma_{i+1}^{-1}...\bar\sigma_g^{-1}\bar a_{g+1}\bar\sigma_g...
\bar\sigma_{i+1}\sigma_{i+2}\sigma_i\sigma_{i+1},\tag5 $$
$$...,$$
$$\bar a_{g-1}=\bar\sigma_{g}^{-1}\bar\sigma_{g-1}^{-1}\bar\sigma_{g+1}^{-1}
\bar\sigma_g^{-1}\bar a_{g+1}\bar\sigma_{g}\bar\sigma_{g+1}\bar\sigma_{g-1}
\bar\sigma_{g},$$
From the equations (4) and (5) we get
$$\bar\sigma_g\bar\sigma_{g+1}...\bar\sigma_{i}\bar\sigma_{i+1}^2\bar\sigma_i
...\bar\sigma_{g+1}\bar\sigma_g\bar a_{g+1}=\bar a_{g+1}\bar\sigma_g
\bar\sigma_{g+1}...\bar\sigma_{i}\bar\sigma_{i+1}^2\bar\sigma_i...
\bar\sigma_{g+1}\bar\sigma_g.\tag6$$
We multiply  on the right side both parts of the equality (6) by
$\bar\sigma_{g+1}$ and then use the braid group relations and the relation
$\bar\sigma _{g+1}\bar a_{g+1}=\bar a_{g+1}\bar\sigma_{g+1}$. We have the
following equalities for the left part of (6):
$$\bar\sigma_{g+1}\bar\sigma_g\bar\sigma_{g+1}...\bar\sigma_{i}
\bar\sigma_{i+1}^2\bar\sigma_i...\bar\sigma_{g+1}\bar\sigma_g\bar a_{g+1}=
\bar\sigma_{g}\bar\sigma_{g+1}\bar\sigma_{g}...\bar\sigma_{i}
\bar\sigma_{i+1}^2\bar\sigma_i...\bar\sigma_{g+1}\bar\sigma_g
\bar a_{g+1}=...$$
$$=\bar\sigma_{g}\bar\sigma_{g+1}\bar\sigma_{g}...\bar\sigma_{i+1}
\bar\sigma_{i}\bar\sigma_{i+1}^2\bar\sigma_i...\bar\sigma_{g+1}\bar\sigma_g
\bar a_{g+1}=$$
$$=\bar\sigma_{g}\bar\sigma_{g+1}\bar\sigma_{g}...\bar\sigma_{i+2}
\bar\sigma_{i}^2\bar\sigma_{i+1}\bar\sigma_i^2\bar\sigma_{i+2}...
\bar\sigma_{g+1}\bar\sigma_g\bar a_{g+1}=$$
$$=\bar\sigma_{g}\bar\sigma_{g+1}\bar\sigma_{g}...\bar\sigma_{i}^2
\bar\sigma_{i+1}\bar\sigma_{i+2}\bar\sigma_{i+1}\bar\sigma_i^2...
\bar\sigma_{g+1}\bar\sigma_g\bar a_{g+1}=...=$$
$$=\bar\sigma_{g}\bar\sigma_{g-1}...\bar\sigma_{i}^2\bar\sigma_{i+1}...
\bar\sigma_g\bar\sigma_{g+1}\bar\sigma_g...\bar\sigma_{i}^2\bar\sigma_{i+1}
...\bar\sigma_g=\tau_i\tau_{i+1}...\tau_g\bar\sigma_{g+1}\tau_i\tau_{i+1}...
\tau_g\bar a_{g+1}.$$
Analogous transformations for the right side of (6) give the equation
$$\tau_{i,g}\bar\sigma_{g+1}\tau_{i,g}\bar a_{g+1}=\bar a_{g+1}\tau_{i,g}
\bar\sigma_{g+1}\tau_{i,g}. \quad \square $$
\enddemo
\proclaim{Corollary 1} The following relations are fulfilled in the monoid
$SB_{n}^g$:
$$\tau_{i,g}\sigma_{1}\tau_{i,g}a_{1}=a_{1}\tau_{i,g}\sigma_{1}\tau_{i,g},
\quad i\leq g. \quad \square$$
\endproclaim
\remark{Remark} The analogous relations are fulfilled in the braid group
$Br_{g+n}$ and hence in the braid group in a handlebody $Br_{n}^g$:
$$\tau_{i,m}\bar\sigma_{m+1}\tau_{i,m}\bar\sigma_{m+1}=
\bar\sigma_{m+1}\tau_{i,m}\bar\sigma_{m+1}, \quad i\leq m.$$ The proof is
straightforward using the induction on $m-i$.
\endremark
There exist maps which are right inverses of the homomorphism of the
canonical inclusion
$$j_k: Br_k\to SG_k.$$
One of them
$$h_k: SG_k\to Br_k,$$
is defined by the formulas:
$$h_k(\bar\sigma_i)=\bar\sigma_i,$$
$$h_k(\bar a_i)=e,$$
and another one
$$h_k^\prime: SG_k\to Br_k,$$
by the following action on generators:
$$h_k^\prime(\bar\sigma_i)=\bar\sigma_i,$$
$$h_k^\prime(\bar a_i)=\bar\sigma_i.$$
The image of the composition
$$SG_n^g\buildrel\roman{Incl}\over\longrightarrow SG_{g+n}
\buildrel h_{g+n}\over\longrightarrow Br_{g+n}$$
lies in $Br_n^g$. We define a homomorphism $h_n^g$ to make the following
diagram commutative
$$\CD
SG_n^g@>\roman{Incl}>> SG_{g+n}\\
\downarrow h_n^g && \downarrow  h_{g+n}\\
Br_n^g@>\roman{Incl}>> Br_{g+n}. \\
\endCD$$
These homomorphisms fit for the following commutative diagram
$$\CD
SG_n@>\roman{Incl}>> SG_n^g@>\roman{Incl}>> SG_{g+n}\\
\downarrow h_n && \downarrow h_n^g && \downarrow  h_{g+n}\\
Br_n@>\roman{Incl}>>Br_n^g@>\roman{Incl}>> Br_{g+n}. \\
\endCD$$
Let $\roman{deg}_k$ be the homomorphism
$$SG_k\to\Bbb Z$$
which assigns to each element of the group $SG_k$ the sum of degrees of the
generators $\bar a_i$, which occur in this element. We define a homomorphism
$$\roman{deg}_n^g: SG_n^g\to\Bbb Z$$
as the following composition
$$SG_n^g\buildrel\roman{Incl}\over\longrightarrow SG_{g+n}
\buildrel\roman{deg}_n^g\over\longrightarrow\Bbb Z.$$
One defines a homomorphism
$$\Bbb Z\to SG_n^g$$
as an inclusion of the cyclic group as a subgroup generated by the element
$a_1$.
Then the composition
$$\Bbb Z\to SG_n^g\buildrel\roman{deg}_n^g\over\longrightarrow\Bbb Z$$
is equal to identity.

\specialhead \S 3. Vassiliev algebra for handlebodies
\endspecialhead

Let $K$ be a commutative ring with unit, $K[\epsilon]$ is a polynomial ring
on one variable $\epsilon$ over $K$, $KSB_n^g$ denotes the semigroup algebra
of the monoid (semigroup with the unit) $SB_n^g$ over $K$.
\par {\bf Definition 1}. We call {\it Vassiliev algebra} $V_n^g$ the
factor-algebra of $KSB_n^g\otimes K[\epsilon]$ by the ideal generated by the
relations
$$\sigma_i-\sigma_i^{-1}=\epsilon a_i.$$
\remark{Remark} Definitions and results of the paper of J.~Baez [5] follow if
we take $g=0$. We omit writing zero in these cases.
\endremark
The homomorphism of inclusion
$$SB_n^g\to SB_{g+n}$$
induces the homomorphism of algebras:
$$\alpha:V_n^g\to V_{n+g},$$
such that $\alpha(\sigma_i)=\bar\sigma_{g+i}$, $\alpha(a_i)=\bar a_{g+i}$,
$\alpha(\tau_k)= \bar\sigma_g\bar\sigma_{g-1}...\bar\sigma_{k+1}
\bar\sigma_{k}^2\bar\sigma_{k+1}^{-1}...\bar\sigma_{g-1}^{-1}
\bar\sigma_g^{-1}.$

We define the homomorphism
$$v^g:KBr_n^g\to V_n^g$$
by the formulas:
$$v^g(\sigma_i)=\sigma_{i},$$
$$v^g(\tau_i)=\tau_i.$$
Let $K(\epsilon)$ be the algebra of Laurent polynomials in $\epsilon$ over
$K$.
\proclaim{Proposition 2} The homomorphism
$$v^g\otimes 1:KBr_n^g\otimes K(\epsilon)\to V_n^g\otimes_{K[\epsilon]}
K(\epsilon)$$
is an isomorphism.
\endproclaim
\demo{Proof} The generators of $V_n^g\otimes_{K[\epsilon]}K(\epsilon)$ are
the images of the elements from $KBr_n^g\otimes K(\epsilon)$: for $\sigma_i$
and $\tau_i$ this follows from definition and for $a_j$ from the formula
$$(v^g\otimes 1)((\sigma_i-\sigma_i^{-1})\otimes\epsilon^{-1})=a_j\otimes 1.
$$
So $v^g\otimes 1$ is an epimorphism. Let us consider the commutative diagram
$$\CD
KBr_n^g\otimes K(\epsilon)@>v^g\otimes 1>>V_n^g\otimes_{K[\epsilon]}
K(\epsilon)\\
\downarrow \roman{Incl}\otimes 1 && \downarrow \alpha\\
KBr_{n+g}\otimes K(\epsilon)@>v\otimes 1>>V_{n+g}\otimes_{K[\epsilon]}
K(\epsilon).\\
\endCD$$
The homomorphism $\roman{Incl}\otimes 1$ is a monomorphism, because all the
modules are free. The homomorphism $v\otimes 1$ is a monomorphism, because it
has evidently defined right inverse. Hence $v^g\otimes 1$ is a monomorphism
and then an isomorphism.
$\square$
\enddemo
\remark{Remark} This Proposition is more important in the case of a
handlebody than in the classical case because the presentation of the
Vassiliev algebra in this case is not so evident.
\endremark
\specialhead \S 4. Invariants of links in handlebodies
\endspecialhead
Let $\Cal L$ be a $K$-valued link invariant by what we mean an ambient
isotopy invariant of oriented links in a handlebody of a genus $g$. It
uniquely extends to a $K(\epsilon)$-valued invariant of generalized links
admitting transverse double points by the standard procedure
$$\Cal L(L_+)-\Cal L(L_-)=\epsilon\Cal L(L_\times),$$
where $L_+$, $L_-$ and $L_\times$ denote link diagrams with a right handed
crossing, left handed crossing, and an intersection respectively, at a given
point, the rest of the diagram being the same. Usually $\epsilon$ is supposed
to be equal to 1 in this equation. An invariant $\Cal L$ is called to be of
{\it degree} $d$ if it vanishes on all generalized links with at least $d+1$
self-intersections. A link invariant of degree $d$ for some $d$ is said to be
of {\it finite type}.

For all $n$ there are inclusions $V_n^g\subset V_{n+1}^g$ and
$KBr_n^g\subset KBr_{n+1}^g$. We consider the direct limits $V_{\infty}^g$
and $KBr_{\infty}^g$ and omit index $\infty$ in notations.
\par {\bf Definition 2}. We call a {\it Markov trace} on $V^g$ to be a
$K[\epsilon]$-linear homomorphism to some $K[\epsilon]$-module $E$
$$T:V^g\to E,$$
satisfying the following conditions
$$T(xy)=T(yx) \ \ \text{for all} \ \ x,y\in V^g,$$
$$T(x\sigma_n^{\pm 1})=zT(x) \ \ \text{for some} \ \ z\in K \ \text{and for
all} \ x\in V^g_n\subset V^g.$$
Similarly is defined the {\it Markov trace} to some $K$-module $M$
$$tr:KBr^g\to M.$$

This definition of Markov traces is analogous to the classical one for usual
braids
$$tr:KBr\to M.$$
\par {\bf Definition 3}. A Markov trace
$$T:V^g\to K[\epsilon]$$
is said to be {\it homogeneous of degree} $d$ if for every $x\in KBr^g$ its
image $T(v(x))\in K[\epsilon]$ is a homogeneous polynomial on $\epsilon$ of
degree $d$.

For a braid $x\in Br^g_n$ we denote by $\hat x$ its closure with respect to
the last $n$ strings. It is considered as a link in a handlebody of the
genus $g$. For a link $L$ we denote by $L\cup\circ$ the unlinked union with
the unknot (which is also supposed to be unlinked with the holes of the
handlebody).
\proclaim {Theorem 1} There is a one-to-one correspondence between $K$-valued
link invariants $\Cal L$ of degree $d$ such that for some $z\in K^*$ (group
of units of the ring $K$)
$$z\Cal L(L\cup\circ)=\Cal L(L)\tag7$$
and Markov traces
$$T:V^g\to K[\epsilon]$$
that are homogeneous of degree $d$. The invariant $\Cal L$ determines the
trace $T$, and conversely by the formula
$$T(v^g(x))=\epsilon^dz^{n-1}\Cal L(\hat x), \ \text{for} \ x\in Br^g_n.\tag8
$$
\endproclaim
\demo{Proof} Let $M$ be $K$-module, $\Cal L$ is a $M$-valued link invariant,
satisfying the equation (7). Then by the relative version of the Markov's
theorem proved by Sofia Lambropoulou and Colin Rourke [11, Theorem 4.7] there
exists a Markov trace
$$tr:KBr^g\to M,$$
given by the formula
$$tr(x)=z^{n-1}\Cal L(\hat x), \ \text{for all} \ x\in Br^g_n.\tag9$$
We note that this trace is well-defined because if $y\in KBr^g_{n+1}$ is the
image of $x\in Br^g_n$ under the canonical inclusion
$KBr_n^g\subset KBr_{n+1}^g$, then
$$tr(y)=z^{n}\Cal L(\hat y)=z^{n}\Cal L(\hat x\cup\circ)=
z^{n-1}\Cal L(\hat x)=tr(x).$$
Moreover, $tr$ is a trace because $\widehat{xy}=\widehat{yx}$ for all
$x, y\in Br^g_n$. We also have
$$tr(x\sigma_n^{\pm 1})=z^{n}\Cal L(\widehat{x\sigma_n^{\pm 1}})=
z^{n}\Cal L(\hat x)=ztr(x).$$
So, equation (9) gives a one-to-one correspondence between $M$-valued link
invariants satisfying (7) and Markov traces
$$tr:KBr^g\to M.$$
Now let $\Cal L$ be $K$-valued link invariant of degree $d$ satisfying
equation (7). Then there is a unique Markov trace
$$tr:KBr^g\to K$$
given by the equation (9). Then by the Proposition 2 there exists a map
$$\tilde T:V^g\otimes_{K[\epsilon]}K(\epsilon)\to K(\epsilon),$$
such that
$$\tilde T(v\otimes 1)=\epsilon^d(tr\otimes 1):KBr_n^g\otimes K(\epsilon)\to
K(\epsilon).$$
Let $T_0$ be defined as the composition
$$V^g\cong V^g\otimes_{K[\epsilon]}K[\epsilon]\buildrel 1\otimes\roman{Incl}
\over\longrightarrow V^g\otimes_{K[\epsilon]}K(\epsilon)\buildrel\tilde T
\over\longrightarrow K(\epsilon).$$
It follows that
$$T_0(v(x))=\epsilon^dtr(x), \ \ x\in KBr^g.$$
If an element $x\in V^g$ can be written as a product of $l$ generators $a_i$
and an arbitrary number of the generators $\sigma _j$ or $\tau_k$, then
$T_0(x)=c\epsilon^{d-l}$ for $c\in K$. Moreover, if $l>d$, then $T_0(x)=0$,
since $\Cal L$ is invariant of degree $d$. Hence
$$T_0(x)=\sum_{i=0}^dc_i\epsilon^i, \ \ c_i\in K,$$
for arbitrary $x\in V^g$. This means that $T_0$ factors through a map
$$T:V^g\to K[\epsilon].$$
It is easy to check that $T$ satisfies necessary conditions.
\par Conversely, suppose that $T:V^g\to K[\epsilon]$ is a homogeneous Markov
trace of degree $d$. We define a Markov trace $tr:KBr^g\to K$ by the formula
$$Tv=\epsilon^dtr.$$
Then equations (8) and (9) define a link invariant. We need to prove that it
is of degree $d$. Let $L$ be a link with $l$ self-intersections, and let
$x_0\in SB^g_n$ be such that its closure $\hat x_0$ is ambient isotopic to
$L$. We use here the generalizations of the theorem of J.~Alexander and the
theorem of A.~A.~Markov for singular links and braids proved by Joan Birman
[6] and Bernd Gemein [10]. Let $x_1$ be
the image of $x_0$ in $V^g$. Then
$$x_1=y_1a_{i_1}y_2a_{i_2}...a_{i_l}y_{l+1},$$
where the elements $y_{i}$ do not have entries of generators $a_j$.
Define $x\in KBr^g$ by the formula
$$x=y_1(\sigma_{i_1}-\sigma_{i_1}^{-1})y_2(\sigma_{i_2}-\sigma_{i_2}^{-1})...
(\sigma_{i_l}-\sigma_{i_l}^{-1})y_{l+1}.$$
Then
$$\epsilon^dtr(x)=Tv(x)=\epsilon^lT(x_1).$$
Since $tr(x)\in K$ and $T(x_1)\in K[\epsilon]$, if $l>d$ we must have
$tr(x)=0$. By construction
$$\Cal L(L)=z^{1-n}tr(x),$$
so it follows that $\Cal L(L)=0$.
$\square$
\enddemo
As an example of the situation considered above let
$$T_S:\Bbb Z Br\to\Bbb Z(\epsilon)$$
be the trace constructed by V.~Turaev [13, Theorem 4.2.1]. It is a Markov
trace with $z=1$:
$$T_S(x\sigma_n^{\pm 1})=T_S(x) \ \text{for all} \ x\in Br_n\subset Br.$$
For $i=0,...,g$; we define a trace
$$T_{S,i}:\Bbb Z Br^g\to\Bbb Z(\epsilon)$$
as the composition
$$\Bbb Z Br^g\buildrel\phi\over\longrightarrow\Bbb Z Br^{g-i}\to\Bbb Z Br
\buildrel T_S\over\longrightarrow\Bbb Z(\epsilon),$$
where the first map is defined by the formulas
$$\phi(\tau_j)=e, \ \text{for} \ j=1,...,i; $$
$$\phi(\tau_j)=\tau_{j-i}, \ \text{for} \ j=i+1,...,g; $$
$$\phi(\sigma_k)=\sigma_k, \ \text{for all} \ k; $$
and the second one is the canonical inclusion. It follows from the
presentations (2) of the groups $Br^g_n$ that the homomorphism $\phi$ is
defined correctly. Then we apply the standard substitution
$q=\exp(\epsilon)$ and obtain the unique Markov trace
$$T_i:V^g\to\Bbb Z[[\epsilon]],$$
such that $T_iv^g=T_{S,i}$. We write
$$T_i=\sum_{d=0}^\infty T_{i,d},$$
where each
$$T_{i,d}:V^g\to\Bbb Z[\epsilon]$$
is a Markov trace homogeneous of degree $d$ with
$$T_{i,d}(x\sigma_n^{\pm 1})=T_{i,d}(x) \ \text{for all} \ x\in V^g_n
\subset V^g.$$
So, each $T_{i,d}$ gives a link invariant $\Cal L_{i,d}$ of degree $d$, such
that
$$\Cal L_{i,d}(\hat x)=\epsilon^{-d}T_{i,d}(v^g(x)), \ \text{for} \ x\in
Br^g_n.$$

\enddocument